\documentclass[]{gNST2e}

\usepackage[latin1]{inputenc}

\newcommand\dontcutmath{\relpenalty=10000\relax\binoppenalty=10000\relax}
\newcommand{\1}{\textnormal{1\hspace{-.55ex}\mbox{l}}}

\begin{document}\dontcutmath

\title{Regenerative block empirical likelihood\\ for Markov chains}
\date{\today }
\author{Hugo Harari-Kermadec
\\\vspace{6pt}
\em{ENS-Cachan \& Universit\'e Paris~1, SAMM}
}
\maketitle

\begin{abstract}
Empirical likelihood is a powerful
semi-parametric method increasingly investigated in the literature. However,
most authors essentially focus on an i.i.d. setting. In the
case of dependent data, the classical empirical likelihood method cannot be
directly applied on the data but rather on blocks of consecutive
data catching the dependence structure.
Generalization of empirical likelihood based on the construction
of blocks of increasing nonrandom length have been proposed for time series satisfying mixing conditions.
 Following some recent
developments in the bootstrap literature, we propose a
generalization for a
large class of Markov chains, based on small blocks of various
lengths. Our approach makes use of the regenerative structure of
Markov chains, which allows us to construct blocks which are almost
independent (independent in the atomic case).  We obtain the
asymptotic validity of the method for positive recurrent Markov
chains and present some simulation results.
\end{abstract}
\textbf{Keywords: Nummelin splitting technique, time series, Empirical Likelihood
\\
MSC codes: 62G05, 62F35, 62F40}
\vfill
{\flushright \small
Hugo Harari-Kermadec\\
Dpt Eco-gestion, ENS-CACHAN\\
61, avenue du Président Wilson\\
94235 Cachan cedex, France\\
\ \\
phone : +331 47 40 20 95\\
fax : +331 47 40 24 60\\
email : hugo.harari@ens-cachan.fr

}

\newpage

\section{Introduction}

Empirical Likelihood (EL), introduced by \cite{ow88}, is a powerful
semi-parametric method. It can be used in a very general setting and
leads to effective estimation, tests and confidence intervals. This
method shares many good properties with the conventional parametric
log-likelihood ratio: both statistics have $\chi^2$ limiting
distribution and are Bartlett correctable, meaning that the error
can be reduced from $\mathcal{O}(n^{-1})$ to $\mathcal{O}(n^{-2})$
by a simple adjustment. An additional property of EL is that the corresponding confidence intervals and
tests do not rely on an estimator of the variance. This last property is specially noticeable for
dependent data, since estimating the variance is then a challenging issue. 

Owen's framework has been intensively studied in the 90's~\citep[see][for an overview]{ow01}, leading to many generalizations
and applications, but mainly for an i.i.d. setting. Some adaptations of EL for dependent data have been introduced, such as the Block Empirical Likelihood (BEL) of~\cite{kitamura97} 
for weakly dependent processes or the subject-wise and elementwise empirical likelihoods of~\cite{wang2010} for longitudinal data.
BEL is inspired by similarities with the bootstrap methodology.
Kitamura proposed to apply the empirical likelihood framework not
directly on the data but on blocks of consecutive data, to catch the
dependence structure. This idea, known as \textit{Block Bootstrap} (BB)
or blocking technique~\citep[in the probabilistic literature, see][for references]{doukhan04} goes back
to~\cite{ku89} in the bootstrap literature and has been intensively
exploited in this field~\citep[see][for a survey]{La03}. However,
the BB performance has been questioned, see
\cite{GK96} and \cite{horowitz03}. Indeed it is known that the
blocking technique distorts the dependence structure of the data
generating process and its performance strongly relies on the choice
of the block size. From a
theoretical point of view, the assumptions used to prove the
validity of BB and BEL are generally strong: it is assumed that the process is
stationary and satisfies some strong-mixing properties (some non-stationary processes can nevertheless be handled, see \cite{synowiecki07} for example). In addition
to having a precise control of the coverage probability of the
confidence intervals, we have to assume that the strong mixing
coefficients are exponentially decreasing \citep[see][]{La03,kitamura97}. Moreover, the choice of
the tuning parameter (the block size) may be quite difficult from
a practical point of view.

In this paper, we focus on generalizing empirical likelihood to
Markov chains. Questioning the restriction implied by the Markovian
setting is a natural issue. It should be mentioned that homogeneous
Markov chain models cover a huge number of time series models. In
particular, a Markov chain can always be written in a nonparametric
way: $X_i=h(X_{i-1},\cdots,X_{i-p},\varepsilon_i)$, where
$(\varepsilon_i)_{i\geq0}$ is i.i.d. with density $f$ and,
for $i>0$, $\varepsilon_i$ is independent of $(X_k)_{0\leq k < i}$~\citep[see][]{kallenberg02}.  Note that both $h$ and $f$ are unknown
functions. Such representations explain why, provided that $p$ is
large enough, any process of length $n$ can be generated by a
Markov chain, see \cite{knight75}. Note also that a Markov chain
may not be necessarily strong-mixing. For instance, the
simple linear model $X_i=\frac12 (X_{i-1}+\varepsilon_i)$ with
$\mathbb P(\varepsilon_i=1)=\mathbb P(\varepsilon_i=0)=\frac12$ is
not strong-mixing~\citep[see][for results on dependence
in Econometrics]{doukhan04}.
\cite{doukhan04} gives many classical econometric models that can be seen as Markovian:
ARMA, ARCH and GARCH processes, bilinear and threshold models.

Our approach is also inspired by some recent developments in the bootstrap literature on Markov chains: instead
of choosing blocks of constant length, we use the Markov chain
structure to choose some adequate cutting times and then we obtain
blocks of various lengths. This construction, introduced
in~\cite{BC04a}, catches the dependence structure. It
is originally based on the existence of an atom for the chain i.e.
an accessible set on which the transition kernel is constant~\citep[see][chapter~1.5]{MT96}. The existence of an atom allows us to cut the chain
into regeneration blocks, separated from each other by a visit to
the atom. These blocks (of random lengths) are independent by the
strong Markov property. Once these blocks are obtained, the
\textit{Regenerative Block-Bootstrap (RBB)} consists in resampling
the data blocks to build new regenerative processes. The rate of convergence of the pivotal statistic
obtained by resampling these blocks
($\mathcal{O}(n^{-1+\varepsilon})$) is better than the one obtained
for the Block Bootstrap ($\mathcal{O}(n^{-3/4})$) and is close to
the classical rate $\mathcal{O}(n^{-1})$ obtained in the i.i.d.
case, see~\cite{GK96} and~\cite{La03}.

These improvements suggest that a version of the empirical
likelihood (EL) method based on such blocks could yield improved
results in comparison to the method presented in~\cite{kitamura97}.
Indeed it is known that EL enjoys somehow the same properties in
terms of  accuracy as the bootstrap but without any Monte-Carlo step.
The main idea is to consider the renewal blocks as independent
observations and to follow the empirical likelihood method. Such
a program is made possible by transforming the original problem based
on moments under the stationary distribution into an equivalent
problem under the distribution of the observable blocks (via Kac's
Theorem). The advantages of the method proposed in this paper are at
least twofold:  first the construction of the blocks is automatic
and entirely determined by the data: it leads to a unique version of
the empirical likelihood program. Second there is not need to ensure
stationarity nor any strong mixing condition to obtain a better
coverage probability for the corresponding confidence regions.

Assuming that the chain is atomic is a strong restriction of this
method. This hypothesis holds for discrete Markov chains
and queuing (or storage) systems returning to a stable state (for
instance the empty queue): see chapter 2.4 of~\cite{MT96}. However
this method can be extended to the more general case of Harris
chains. Indeed, any chain having some recurrent properties can be
extended to a chain possessing an atom which then enjoys some
regenerative properties. Nummelin gives an explicit construction of
such an extension that we recall in Section~\ref{harris}
\citep[see][]{Nu78,AN78}. In \cite{BC04c}, an extension
of the RBB procedure to general Harris chains based on the Nummelin'
splitting technique is proposed (\textit{the Approximate
Regenerative Block-Bootstrap, ARBB}). One purpose of this paper is
to prove that these approximatively regenerative blocks can also be
used in the framework of empirical likelihood and lead to consistent
results.

The outline of the paper is the following. In
Section~\ref{preliminary}, notations are set out and key concepts of
the Markov atomic chain theory are recalled. In
Section~\ref{regenerative}, we present how to construct regenerative
data blocks and confidence regions based on these blocks. We give the main
properties of the corresponding asymptotic statistics.
In Section~\ref{harris} the Nummelin splitting technique is shortly
recalled and a framework to adapt the regenerative empirical
likelihood method to general Harris chains is proposed. We
essentially obtain consistent results but also briefly discuss
test and higher order properties. In Section~\ref{illu}, we present some
moderate sample size simulations.

\section{Preliminary statement}\label{preliminary}

\subsection{Framework}

For the sake of simplicity, we use the same notations as~\cite{BC04d} when possible. We consider a chain $X=\left(X_{i}\right)_{i\in\mathbb{N}}$ on a state space $(E,\mathcal E)$, with initial distribution $\nu$ and transition probability $\Pi$. For a set $B\in\mathcal{E}$ and $i\in\mathbb{N}$, we thus denote
$$
X_{0}\sim\nu\text{ and }\mathbb{P}(X_{i}\in B\mid X_{0},...,\
X_{i-1})=\Pi(X_{i-1},\;B)\text{ a.s. .}
$$

Recurrence properties will be important in the following.
An irreducible chain is said \textit{positive recurrent} when it admits an invariant probability :
$$
\exists \mu \text{ probability measure on }E,\ 
\mu\Pi=\mu,\text{ where } \mu\Pi(\cdot)=\int_{x\in E}\mu(dx)\Pi\left( x,\cdot\right).
$$

We assume that the chain is aperiodic (i.e. $X$ is not cyclic) and that there exists
a measure $\psi$ such as the chain is $\psi$-irreducible.
This simply means that for any starting state $x$ in $\mathcal E$ and any set $A$ of positive $\psi$-measure, the chain visits $A$ with probability~1.
A $\psi$-irreducible chain is said \textit{Harris recurrent}
if every measurable set with positive $\psi$-measure visited once is visited infinitely often with probability 1.

In what follows, $\mathbb{P}_{\nu}$ and $\mathbb{P}_{x}$ (for $x$
in $E$) denote respectively the probability measure when $X_{0}\sim\nu$ and $X_{0}=x$. The indicator function of an event $\mathcal A$ is denoted by $\1_{\mathcal{A}}$. The corresponding expectations are denoted $\mathbb{E}_{\nu}(\cdot)$, $\mathbb{E}_{x}[\cdot]$ and $\mathbb{E}_{\mathcal{A}}[\cdot]$.
For further details and traditional properties of Markov chains, we refer to \cite{Re84} or \cite{MT96}. 

Notice that the chain $X$ is not supposed stationary
(since $\nu$ may differ from $\mu$) nor strong-mixing. To simplify the exposition, we do not treat in this paper the fully non-stationary case corresponding to null recurrence.
Results in that direction may be found in~\cite{tjostheim90}.

\subsection{Atomic Markov chains}

Assume that the chain is $\psi$-irreducible and possesses an
accessible atom, i.e. a set $A$ with $\psi(A)>0$ such
that $\Pi(x,.)=\Pi(y,.)$ for all $x,y$ in $A$. The class of atomic Markov
chains contains not only chains defined on a countable state space
but also many specific Markov models used to study queuing systems
and stock models \citep[see][for models involved in queuing theory]{As87}. In the discrete case, any recurrent state is an accessible
atom: the choice of the atom is thus left to the statistician who
can for instance use the most visited point. In many other
situations the atom is determined by the structure of the model (for
a random walk on $\mathbb{R}^+$, with continuous increment, 0 is the
only possible atom).

Denote by $\tau_{A}=\tau_{A}(1)=\inf\left\{ k\geq1,\ X_{k}\in
A\right\}$ the hitting time of the atom $A$ (the first visit) and,
for $j\geq2$, denote by $\tau_{A}(j)=\inf\left\{ k>\tau_{A}(j-1),\
X_{k}\in A\right\} $ the successive return times to $A$. The
sequence $\left( \tau_{A}(j)\right) _{j\geq1}$ defines the
successive times at which the chain forgets its past, called
\textit{regeneration times}. Indeed, the transition probability
being constant on the atom, $X_{\tau_{A}+1}$ only depends on the
information that $X_{\tau_{A}}$ is in $A$ and not any more on the
actual value of $X_{\tau_{A}}$ itself.

For any initial distribution $\nu$, the sample path of the chain
may be divided into blocks of random length corresponding to
consecutive visits to $A$:
\[
B_{j}=(X_{\tau_{A}(j)+1},...,\ X_{\tau_{A}(j+1)}).
\]
The sequence of blocks $(B_{j})_{1\leq j< \infty}$ is then i.i.d.\ by the strong Markov property \citep[][page~73]{MT96}.
Notice that the block $B_0=(X_1,...,\ X_{\tau_A})$ is independent of the other blocks, but does not have the same
distribution, because it depends on the initial distribution $\nu$.

Let $m:E\times\mathbb R^p \rightarrow\mathbb{R}^r$ be a measurable function and
$\theta_0$ be the true value of some parameter
$\theta\in\mathbb{R}^p$ of the chain, given by an estimating
equation on the invariant measure $\mu$:
\begin{equation}
\mathbb{E}_\mu[m(X,\theta_0)]=0.\label{atomic_moment}
\end{equation}
Dimensions are of importance: the number of constraints $r$ must be at least equal to the number of parameters $p$, for identification reasons.
The estimation of the mean is a just-identified case ($r=p$):
 $\theta_0=\mathbb{E}_\mu[X]$ and $m(X,\theta)=X-\theta$.

In this framework, Kac's Theorem, stated below \citep[][Theorem~10.2.2]{MT96}
allows us to write functionals of the stationary
distribution $\mu$ as functionals of the distribution of a
regenerative block.
\begin{theorem}[Kac]
Let $X$ be an aperiodic, $\psi$-irreducible Markov chain with an accessible atom $A$.
$X$ is positive recurrent if and only if $\mathbb{E}_{A}[\tau _{A}]<\infty $.
In such a case, $X$ admits an unique invariant probability distribution $\mu$, the Pitman's occupation measure given by
\[
\mu(F)=\mathbb{E}_{A}\left[\sum_{i=1}^{\tau _{A}}{\emph\1}_{X_{i}\in F}\right]/
\mathbb{E}_{A}[\tau _{A}],\text{ for all }F\in \mathcal{E}.
\]
\end{theorem}

In the following we denote
$$
M(B_j,\theta)=\sum_{i=\tau_A(j)+1}^{\tau_A(j+1)} m(X_i,\theta)$$ so
that we can rewrite the estimating equation~(\ref{atomic_moment})
as:
\begin{equation}\label{momentM}
\mathbb{E}_A[M(B_j,\theta_0)]=0.
\end{equation}

Kac's Theorem allows us to use the
decomposition of the chain into independent blocks to
obtain limit theorems for atomic chains. See for example
\cite{MT96} for the Law of Large Numbers (LLN, page~415), Central Limit
Theorem (CLT, page~416), Law of Iterated Logarithm (page~416), \cite{Bol82} for the
Berry-Esseen Theorem and \cite{BC04a} for Edgeworth expansions. These
results are established under hypotheses related to the
distribution of the $B_{j}$'s :

\textbf{Return time conditions:}
\begin{align*}
\textbf{H0}(\kappa) & :\ \mathbb{E}_{A}[\tau_{A}^{\kappa}]<\infty,
\\
\textbf{H0}(\kappa,\ \nu) & :\ \mathbb{E}_{\nu}[\tau_{A}^{\kappa}]<\infty,
\end{align*}
where $\kappa>0$ and $\nu$ is the initial distribution of the chain.
When the chain is stationary and strong mixing, these hypotheses can be related to
the rate of decay of $\alpha$-mixing coefficients $\alpha(p)$,
see~\cite{Bol82}. In particular, the hypotheses are satisfied if
$\sum_{j\geq1} j^\kappa\alpha(j)<\infty$.

\textit{Block-moment conditions:}
\begin{align*}
\textbf{H1}(\kappa,\ m) &
 :\ \mathbb{E}_{A}\left[\left(\sum_{i=1}^{\tau_{A}}\left|\left| m(X_{i},\theta_0)\right|\right|\right)^{\kappa}\right]<\infty, \\
\textbf{H1}(\kappa,\ \nu,\ m) &
 :\ \mathbb{E}_{\nu}\left[\left(\sum_{i=1}^{\tau_{A}}\left|\left| m(X_{i},\theta_0)\right|\right|\right)^{\kappa}\right]<\infty.
\end{align*}

The assumptions on $\nu$ allow to control the first block $B_{0}$.
Equivalence of these assumptions with easily checkable drift
conditions may be found in \cite{MT96}, Appendix~A.

\section{The regenerative case}\label{regenerative}

\subsection{Regenerative Block Empirical Likelihood algorithm}

Let $X_{1},\cdots,\ X_{n}$ be an observation of the chain $X$.
If we assume that we know an atom $A$ for the
chain, the construction of the regenerative blocks is then trivial.
Consider the empirical distribution of the blocks:
$$
\mathbb P_{l_n}=\frac 1{l_n}\sum_{j=1}^{l_n}\delta_{B_j},
$$
where $l_n$ is the number of complete regenerative blocks, and the multinomial distributions
$$
\mathbb Q=\sum_{j=1}^{l_n}q_j\delta_{B_j}, \text{ with } 0<q_j \text{ and } \sum_{j=1}^{l_n}q_j=1,
$$
dominated by $\mathbb P_{l_n}$. To obtain a confidence region, we
will apply \cite{ow90}'s method to the blocks $B_j$: we are going to
minimize the Kullback discrepancy between $\mathbb Q$ and $\mathbb
P_{l_n}$ under the condition~(\ref{momentM}). More precisely, the
\textit{Regenerative Block Empirical Likelihood} is defined in the next 4 steps:

\begin{algorithm}[ReBEL - Regenerative Block Empirical Likelihood construction ]
\ \label{algo1}
\begin{enumerate}
    \item  Count the number of visits to $A$ up to time $n$:
    $l_{n}+1=\sum_{i=1}^{n}{\emph\1}_{X_{i}\in A}$.
    \item  Divide the observed trajectory $X^{(n)}=(X_{1},....,X_{n})$
        into $l_{n}+2$ blocks corresponding to the pieces of the
        sample path between consecutive visits to the atom $A$,
        \begin{align*}
                B_{0}& =(X_{1},...,\text{ }X_{\tau _{A}(1)}),
                \ B_{1}=(X_{\tau _{A}(1)+1},...,\text{ }X_{\tau _{A}(2)}),...,\ \\
                B_{l_{n}}& =(X_{\tau _{A}(l_{n})+1},...,
                \ X_{\tau_{A}(l_{n}+1)}),\ B_{l_{n}+1}^{(n)}
                =(X_{\tau _{A}(l_{n}+1)+1},...,\ X_{n}),
        \end{align*}
        with the convention $B_{l_{n}+1}^{(n)}=\emptyset$
        when $\tau _{A}(l_{n}+1)=n$.
    \item  Drop the first block $B_{0}$ and the last one $B_{l_{n}+1}^{(n)}$
     (possibly empty when~$\tau _{A}(l_{n}+1)=n$).
    \item Evaluate the empirical log-likelihood ratio $r_{n}(\theta )$ (practically on a grid of the set of interest):
        $$
        r_{n}(\theta ) =-\sup_{(q_1,\cdots,q_{l_n})} \left\{ \left. \log\left[\prod_{j=1}^{l_n}l_n q_j\right]\right|
        \sum_{j=1}^{l_n}q_j\cdot M(B_j,\theta)=0,\ \sum_{j=1}^{l_n}q_j=1 \right\}.
        $$
        Using Lagrange arguments, this can be more easily calculated as
        $$
        r_{n}(\theta)=\sup_{\lambda\in\mathbb{R}^p }
        \left\{ \sum_{j=1}^{l_n}\log \left[1+\lambda'M(B_j,\theta)\right] \right\}.
        $$
\end{enumerate}
\end{algorithm}

\begin{remark}[\underline{Small samples} ]
Possibly, if the chain does not visit $A$, $l_n=-1$. Of course
the algorithm cannot be implemented and no confidence interval
can be built. Actually, even when $l_n\geq0$, the algorithm can
be meaningless and at least a reasonable number of blocks are
needed to build a confidence interval. In the positive recurrent
case, it is known that $l_n\sim n/\mathbb{E}_A[\tau_A]$ a.s. and
the length of each block has expectation $\mathbb{E}_A[\tau_A]$.
Many regenerations of the chain should then be observed as soon as
$n$ is significantly larger than $\mathbb{E}_A[\tau_{A}]$. Of
course, the next results are asymptotic, for finite sample
consideration on empirical likelihood methods (in the i.i.d.
setting), refer to \cite{ber-gaut-hara08}.
\end{remark}

The next theorem states the asymptotic validity of ReBEL
in the case $r=p$ (just-identified case). For this, we
introduce the ReBEL confidence region defined as follows:
        $$
        C_{n,\alpha}=\left\{ \theta \in \mathbb{R}^{p}\left|
        \;2\cdot r_{n}(\theta)\leq F^{-1}_{\chi_p^2}(1-\alpha )\right.\right\},
        $$
where $F_{\chi_p^2}$ is the distribution function of a $\chi^2$ distribution with $p$ degrees of freedom.
\begin{remark}[\underline{Scaling factor} ]
Block empirical likelihood methods usually need a scaling factor to compensate the eventual overlap of the blocks, denoted $A_N^{-1}$ in~\cite{kitamura97}, page~2089. The regenerative perspective used here forbid any overlap and therefore avoid such a factor.
\end{remark}

\begin{theorem}\label{th_rebel}
Let $\mu$ be the invariant measure of the chain, let $\theta_0\in \mathbb{R}^p$ be the parameter of interest,
satisfying $\mathbb{E}_\mu[m(X,\theta_0)]=0$. Assume
$$
\Sigma=\mathbb E_A[\tau_A]^{-1}\mathbb{E}_A[M(B,\theta_0)M(B,\theta_0)']
$$
is of full-rank.
If \textbf{H0}$(1,\nu)$, \textbf{H0}$(2)$ and \textbf{H1}$(2,m)$ hold, then
$$2r_{n}(\theta_0)\xrightarrow[n\to\infty]{\mathcal L} \chi_{p}^{2}$$
and therefore
\[
\mathbb{P}_\nu\left(\theta_0\in C_{n,\alpha}\right)\xrightarrow[n\to\infty]{} 1-\alpha.
\]
\end{theorem}

The proof relies on the same arguments as the one for empirical
likelihood based on i.i.d. data. This can be easily understood: our
data, the regenerative blocks, are i.i.d. \citep{ow90,ow01}. The only difference with the classical use of
empirical likelihood is that the length of the data (i.e. the number
of blocks) is a random value $l_n$. However, we have that $n/l_n \to
\mathbb{E}_A(\tau_A)$ a.s. \citep[][page~425]{MT96}. The proof is given in
the appendix.

\begin{remark}[\underline{Convergence rate} ]
Let's make some very brief discussion on the rate of convergence of
this method. \cite{BC04a} shows that the Edgeworth expansion of the
mean standardized by the empirical variance holds up to
$\mathcal{O}_{\nu}(n^{-1})$ (in contrast to what is
expected when considering a variance built on fixed length blocks).
It follows from their result that
$$
\mathbb{P}_\nu\left(2 r_{n}(\theta_0)\leq
u\right)=F^{-1}_{\chi^2_p}(u)+\mathcal{O}_{\nu}(n^{-1})
$$
This is already (without Bartlett correction) better than the
Bartlett corrected empirical likelihood when fixed length
blocks are used \citep{kitamura97}. Actually, we expect, in this atomic
framework, that a Bartlett correction would lead to the same result
as in the i.i.d. case: $\mathcal{O}(n^{-2})$. However, to prove this
conjecture, we should establish an Edgeworth expansion for the
likelihood ratio (which can be derived from the Edgeworth expansion for
self-normalized sums) up to order $\mathcal{O}(n^{-2})$ which is a very
technical task. This is left for further work.
\end{remark}

\begin{remark}[\underline{Change of discrepancy} ]\label{change-div}
Empirical likelihood can be seen as a contrast method based on the
Kullback discrepancy. To replace the Kullback discrepancy by some
other discrepancy is an interesting problem which has led to some
recent works in the i.i.d. case. \cite{newey-smith04} generalized
empirical likelihood to the family of Cressie-Read discrepancies
\citep[see also][]{gug-smith05}. The resulting methodology,
\textit{Generalized Empirical Likelihood}, is included in the
empirical $\varphi$-discrepancy method introduced by
\cite{ber-gaut-hara08} \citep[see also][]{ber-hara-rav07,kitamura06}.

In the dependent case, it should be mentioned that the constant
length blocks procedure has been studied in the case of empirical
Euclidean likelihood by~\cite{lin01}. A method based on the Cressie-Read discrepancies for tilting time series data
has been introduced by~\cite{hall-yao03}. Our proposal, stated here for the Kullback discrepancy only, is straightforwardly compatible with these generalizations (Cressie-Read and $\varphi$-discrepancy).
\end{remark}

An important issue is the behavior of the empirical log-likelihood ratio under a local alternative, i.e. if the moment equation~(\ref{atomic_moment}) is misspecified :
$\mathbb{E}_\mu[m(X,\theta_0)]=\delta/\sqrt n$.
The result states as follows.

\begin{theorem}\label{th_power}
Let $\mu$ be the invariant measure of the chain, let $\theta_0\in \mathbb{R}^p$ be the parameter of interest,
satisfying $\mathbb{E}_\mu[m(X,\theta_0)]=\delta/\sqrt n$. Assume that
$\Sigma$ is of full-rank.
If \textbf{H0}$(1,\nu)$, \textbf{H0}$(2)$ and \textbf{H1}$(2,m)$ hold, then
the empirical log-likelihood ratio has an asymptotic noncentral chi-square distribution with
$p$ degrees of freedom and noncentrality parameter $\delta'\Sigma^{-1}\delta$
$$
2r_{n}(\theta_0)\xrightarrow[n\to\infty]{\mathcal L} \chi^{'2}_p(\delta'\Sigma^{-1}\delta).
$$
\end{theorem}
The proof is postponed to the appendix. It is a classical result that the log-likelihood ratio is asymptotically noncentral chi-square and that the critical order is $n^{-1/2}$. The interesting quantity to study the efficiency
of the method in this context is the noncentrality parameter. \cite{newey85} gives the asymptotic distribution of
the pivotal statistic based on optimally weighted GMM which is a standard tool for dependent data. Unfortunately,
Newey's results are stated in a parametric context and it is therefore impossible to compare them with Theorem~\ref{th_power}.

Nevertheless, ReBEL can easily be compared with the Continuously updated GMM (CUE-GMM)
which is very close to the optimally weighted GMM. CUE-GMM estimators have been shown to coincide
with empirical Euclidean likelihood (EEL), see~\cite{bonnal-renault07}.
The difference between the EL and EEL being just a change of discrepancy (see Note~\ref{change-div}), it is then straightforward to
adapt the proof of Theorem~\ref{th_power} to the case of the EEL. The developments of the pivotal statistics coincide for the two first order and therefore they lead to the same asymptotic distribution in the case of misspecification.
EL is thus as efficient as the optimally weighted GMM.

\subsection{Estimation and the over-identified case}

The properties of empirical likelihood proved by \cite{qin94} can
be extended to our Markovian setting. In order to state the
corresponding results respectively on estimation, confidence region
under over-identification ($r> p$) and hypotheses testing, we
introduce the following additional assumptions. Assume that there
exists a neighborhood $V$ of $\theta_0$ and a real positive function $N$ with
$\mathbb{E}_{\mu}\left[N(X)\right]<\infty$, such that:
\begin{itemize}
\item[\textbf{H2}(a)] $\partial m (x,\theta)/\partial\theta$ is continuous in $\theta$
and bounded in norm by $N(x)$ for $\theta$ in $V$,
\item[\textbf{H2}(b)] $D=\mathbb E_\mu[\partial m (X,\theta_0)/\partial\theta]$ is of full rank,
\item[\textbf{H2}(c)] $\partial^2 m (x,\theta)/\partial\theta\partial\theta'$ is continuous in $\theta$
and bounded in norm by $N(x)$ for $\theta$ in $V$,
\item[\textbf{H2}(d)] $\|m(x,\theta)\|^3$ is bounded by $N(x)$ on $V$.
\end{itemize}

Notice that \textbf{H2}(d) implies in particular the block moment
condition $\textbf{H1}(3,m)$ since by Kac's Theorem
$$\mathbb E_{\mu}\left[\|m(X,\theta)\|^3\right]=\frac{\mathbb
E_A\left[\sum_{i=1}^{\tau_A}\|m(X_i,\theta)\|^3\right]}{\mathbb E_A[\tau_A]}
\leq
\frac{\mathbb
E_A\left[\sum_{i=1}^{\tau_A}N(X_i)\right]}{\mathbb E_A[\tau_A]}
=\mathbb E_{\mu}\left[N(X)\right]<\infty.
$$

Empirical likelihood provides a natural way to estimate $\theta_0$
in the i.i.d. case \citep{qin94}. This can be
straightforwardly extended to Markov chains. The estimator is the
maximum empirical likelihood estimator defined by
$$
\tilde\theta_n=\arg\inf_{\theta\in\Theta} \{r_n(\theta)\}.
$$
The next theorem shows that, under natural assumptions on $m$ and
$\mu$, $\tilde\theta_n$ is an asymptotically Gaussian estimator of
$\theta_0$.
\begin{theorem}\label{atomic_estimator}
Assume that the hypotheses of Theorem~\ref{th_rebel} holds. Under
the additional assumptions~\textbf{H2}(a), \textbf{H2}(b) and
\textbf{H2}(d), $\tilde\theta_n$ is a consistent estimator of
$\theta_0$. If in addition~\textbf{H2}(c) holds, then
$\tilde\theta_n$ is asymptotically Gaussian:

$$
\sqrt{n}(\tilde\theta_n-\theta_0)\xrightarrow[n\to\infty]{\mathcal
L} \mathcal N\left(0,\left(D'\Sigma^{-1}D\right)^{-1}\right).
$$
\end{theorem}

Notice that both $D$ and $\Sigma$ can
be easily estimated by empirical sums over the blocks.
The corresponding estimator for $\left(D'\Sigma^{-1}D\right)^{-1}$ is straightforwardly
convergent by the LLN for Markov chains.

\begin{remark}[\underline{Asymptotic covariance matrix} ]
Our asymptotic covariance matrix $\left(D'\Sigma^{-1}D\right)^{-1}$
is to be compared with the asymptotic covariance matrix
$V_\theta$ of~\cite{kitamura97}'s estimator, which coincide with the
asymptotic covariance matrix of the optimally weighted GMM estimator.
Both matrix are very similar: $V_\theta=(D'\mathcal S^{-1}D)^{-1}$, where $\mathcal S$ is
the  counterpart of our $\Sigma$ for weakly dependent processes:
$$
\mathcal S=\lim_{n\to\infty}n^{-1}\left(\sum_{i=1}^n m(X_i,\theta_0)\right)\left(\sum_{i=1}^n m(X_i,\theta_0)\right)'.
$$
For a process being both weakly dependent and Markovian (and in particular in the i.i.d. case),
$\mathcal S=\Sigma$ and therefore $V_\theta=\left(D'\Sigma^{-1}D\right)^{-1}$.
\end{remark}

The case of over-identification ($r>p$) is an important feature,
specially for econometric applications. In such a case, the
statistic $2r_n(\tilde\theta_n)$ may be considered to test the
moment equation~(\ref{atomic_moment}):
\begin{theorem}\label{atomic_test}
Under the assumptions of Theorem~\ref{atomic_estimator}, if the
moment equation~(\ref{atomic_moment}) holds, then we have
$$2r_n(\tilde\theta_n)\xrightarrow[n\to\infty]{\mathcal L}\chi^2_{r-p}.$$
\end{theorem}

We now turn  to a theorem equivalent to Theorem~\ref{th_rebel}. In the
over-identified case, the likelihood ratio statistic used to
test $\theta=\theta_0$ must be corrected. We now define
$$
W_{1,n}(\theta)=2r_n(\theta)-2r_n(\tilde\theta_n).
$$
The ReBEL confidence region of nominal level $1-\alpha$ in the
over-identified case is now given by
$$
C^1_{n,\alpha}=\left\{ \theta \in
\mathbb{R}^{p}\left|W_{1,n}(\theta) \leq F^{-1}_{\chi_p^2}(1-\alpha
)\right.\right\}.
$$
\begin{theorem}\label{atomic_overid}
Under the assumptions of Theorem~\ref{atomic_estimator},
the likelihood ratio statistic for $\theta=\theta_0$ is asymptotically $\chi_{p}^{2}$:
$$
W_{1,n}(\theta_0)\xrightarrow[n\to\infty]{\mathcal L} \chi_{p}^{2}
$$
and $C^1_{n,\alpha}$ is then an asymptotic confidence region of
nominal level $1-\alpha$.
\end{theorem}

To test a sub-vector of the parameter, we can also build the
corresponding empirical likelihood ratio \citep{qin94,kitamura97,kitamura04,gug-smith05}. Let
$\theta=(\gamma,\beta)'$ be in $\mathbb R^q\times\mathbb
R^{p-q}$, where $\gamma\in\mathbb R^q$ is the parameter of
interest and $\beta\in\mathbb R^{p-q}$ is a nuisance parameter.
Assume that the true value of the parameter of interest is
$\gamma_0$. The empirical likelihood ratio statistic in this case
becomes
$$
W_{2,n}(\gamma)=
2\cdot\left(\inf_{\beta}r_n((\gamma,\beta)'
)-\inf_{\theta}r_n(\theta)\right)
=2\cdot\left(\inf_{\beta}r_n((\gamma,\beta)'
)-r_n(\tilde\theta_n)\right),
$$
and the empirical likelihood confidence region is given by
$$
C^2_{n,\alpha}=\left\{ \gamma \in
\mathbb{R}^{q}\left|W_{2,n}(\gamma) \leq F^{-1}_{\chi_q^2}(1-\alpha
)\right.\right\}.
$$
\begin{theorem}\label{atomic_subvector}
Under the assumptions of Theorem~\ref{atomic_estimator},
$$W_{2,n}(\gamma_0)\xrightarrow[n\to\infty]{\mathcal L}\chi^2_{q}$$
and $C^2_{n,\alpha}$ is then an asymptotic confidence region of
nominal level $1-\alpha$.
\end{theorem}

\section{The case of general Harris chains}\label{harris}

\subsection{Algorithm}\label{algo}
As explained in the introduction, the splitting technique introduced in~\cite{Nu78}
 allows us to extend our algorithm to
general Harris recurrent chains. The idea is to extend the
original chain to a ``virtual" chain with an atom. The splitting technique relies
on the crucial notion of \textit{small set}. Additional definitions are needed: 
a set $S\in\mathcal{E}$ is said to be \textit{small} if
there exist $\delta>0$, a positive integer $q$ and a probability measure $\Phi$ supported by $S$ such that,
for all $x\in S,$ $A\in\mathcal{E}$,
\begin{equation}\label{minor}
\Pi^{q}(x,A)\geq\delta\Phi(A),
\end{equation}
$\Pi^{q}$ being the $q$-th iterate of the transition probability $\Pi$.
Note that an accessible small set
always exists for $\psi$-irreducible chains \citep{JJ67}.

In the case $q>1$, a first step is typically  to reduce the order to 1 by stacking lagged values (an example in given in section~\ref{tgarch}). Nevertheless, this complicates the exposition and the demonstrations since the resulting transition probability has no density and since the splitting technique leads to 1-dependence instead of independence. See~\cite{BCT09} and~\cite{Adamczak08} on that issues.
For simplicity, we assume in the following that $q=1$ and that $\Phi$ has a density $\phi$
with respect to some reference measure $\lambda(\cdot)$.
 
The idea to construct the split chain $\widetilde{X}=(X,W)$ is the
following:
\begin{itemize}
\item if $X_i\notin S$, generate (conditionally to $X_i$) $W_i$ as a Bernoulli random value,
 with probability $\delta$.
\item if $X_i\in S$, generate (conditionally to $X_i$) $W_i$ as a Bernoulli random value,
 with probability $\delta\phi(X_{i+1})/p(X_{i},X_{i+1})$,
\end{itemize}
where $p$ is the transition density of the chain $X$ with respect to $\lambda$. This
construction essentially relies on the fact that under the
minorization condition~(\ref{minor}), $\Pi(x,A)$ may be written on
$S$ as a mixture:
$\Pi(x,A)=(1-\delta)(\Pi(x,A)-\delta\Phi(A))/(1-\delta)+\delta\Phi(A)$,
which is constant (independent of the starting point $x$) when the second component is picked \citep[see][for details]{BC05}.

When constructed this way, the split chain is an atomic Markov
chain, with marginal distribution equal to the original distribution
of $X$ \citep[see][page~427]{MT96}. The atom is then $A=S\times \{1\}$. In
practice, we will only need to know when the split chain visits the
atom, i.e. we only need to simulate $W_i$ when $X_i\in S$.

Those visits to the atom are therefore the date of regeneration of the chain, and the number of visits acts as a sample size. In practice, the choice of the small set is then decisive for the performance of the algorithm. 
A balance needs to be achieved: if $S$ were chosen too large, it would be visited very often, but the minorization condition~(\ref{minor}) would likely be poor and therefore $\delta$ would be small. This would lead to many realization $W_i=0$ and few $W_i=1$. Most of the visits of $X_i$ to the small set would then be wasted since they would not give a regeneration time. This balance is not a curse: it gives a natural data-driven tuning of the small set and prevent from the difficulties rising in the choice of kernel bandwidth for example. For a discussion on the practical choice of the small set, see~\cite{BC05}.

The return time conditions are now defined as uniform moment
condition over the small set:
\begin{align*}
\textbf{H0}(S,\ \kappa) & :\ \sup_{x\in S}\mathbb{E}_{x}[\tau_{S}^{\kappa}]<\infty,
\\
\textbf{H0}(S,\ \kappa,\ \nu) & :\
\mathbb{E}_{\nu}[\tau_{S}^{\kappa}]<\infty.
\end{align*}
The Block-moment conditions become
\begin{align*}
\textbf{H1}(S,\ \kappa,\ m) &
 :\ \sup_{x\in S}\mathbb{E}_{x}\left[\left(\sum_{i=1}^{\tau_{S}}\|m(X_{i},\theta_0)\|\right)^{\kappa}\right]<\infty, \\
\textbf{H1}(S,\ \kappa,\ \nu,\ m) &
 :\ \mathbb{E}_{\nu}\left[\left(\sum_{i=1}^{\tau_{S}}\|m(X_{i},\theta_0)\|\right)^{\kappa}\right]<\infty.
\end{align*}

Unfortunately, the Nummelin technique involves the transition
density of the chain, which is of course unknown in a nonparametric
approach. An approximation $p_n$ of this density can however be
computed easily by using standard kernel methods. This leads us to
the following version of the empirical likelihood program.

\begin{algorithm}[Approximate regenerative block EL construction ]
\ \label{algo2}
\begin{enumerate}
  \item Find an estimator $p_n$ of the transition density (for instance a Nadaraya-Watson estimator).
  \item Choose a small set $S$ and a density $\phi$ on $S$ and evaluate
  $\delta=\min_{x,y\in S}\left\{\frac{p_n(x,y)}{\phi(y)}\right\}$.
  \item When $X$ visits $S$, generate $\widehat W_i$ as a Bernoulli with parameter $\delta\phi(X_{i+1})/p_n(X_{i},X_{i+1})$.
          If $\widehat W_i=1$, the approximate split chain $(X_i,\widehat W_i)=\widehat{X}_{i}$ visits the
          atom $A=S\times\{1\}$ and $i$ is an approximate regenerative
          time. These times define the approximate return times $\widehat\tau_A(j)$.
  \item Count the number of visits to $A$ up to time $n$:
            $\hat{l}_{n}+1=\sum_{i=1}^{n}{\emph\1}_{\widehat{X}_{i}\in A}$.
  \item Divide the observed trajectory $X^{(n)}=(X_{1},....,X_{n})$
          into $\hat{l}_{n}+2$ blocks corresponding to the pieces of the
        sample path between approximate return times to the atom $A$,
        \begin{align*}
                \widehat{B}_{0}& =(X_{1},...,\ X_{\widehat\tau_A(1)}),
                \ \widehat{B}_{1}=(X_{\widehat\tau_A(1)+1},...,\ X_{\widehat\tau_A(2)}),...,\ \\
                \widehat{B}_{\hat{l}_{n}}& =(X_{\widehat\tau_A(\hat{l}_{n})+1},...,
                \ X_{\widehat\tau_A(\hat{l}_{n}+1)}),\ \widehat{B}_{\hat{l}_{n}+1}^{(n)}
                =(X_{\widehat\tau_A(\hat{l}_{n}+1)+1},...,\ X_{n}),
                \end{align*}
                with the convention $\widehat{B}_{\hat{l}_{n}+1}^{(n)}=\emptyset$
                when $\widehat\tau_A(\hat{l}_{n}+1)=n$.
    \item Drop the first block $\widehat{B}_{0},$ and
                the last one $\widehat{B}_{\hat{l}_n+1}^{(n)}$ (possibly empty
                when $\widehat\tau_A(\hat{l}_{n}+1)=n$).
  \item Define
  $$
  M(\widehat B_j,\theta)=\sum_{i=\widehat\tau_A(j)+1}^{\widehat\tau_A(j+1)} m(X_i,\theta).
  $$ Evaluate the empirical log-likelihood ratio $\hat r_{n}(\theta )$ (practically on a grid of the set of interest):
        $$
        \hat r_{n}(\theta ) =-\sup_{(q_1,\cdots,q_{\hat l_n})} \left\{ \left. \log\left[\prod_{j=1}^{\hat l_n}\hat l_n q_j\right]\right|
        \sum_{j=1}^{\hat l_n}q_j\cdot M(\widehat B_j,\theta)=0,\ \sum_{j=1}^{\hat l_n}q_j=1 \right\}.
        $$
        Using Lagrange arguments, this can be more easily calculated as
        $$
        \hat r_{n}(\theta)=\sup_{\lambda\in\mathbb{R}^p }
        \left\{ \sum_{j=1}^{\hat l_n}\log \left[1+\lambda'M(\widehat B_j,\theta)\right] \right\}.
        $$
\end{enumerate}
\end{algorithm}

\subsection{Main theorem}
The practical use of this algorithm crucially relies on
the preliminary computation of a consistent estimator of the
transition density. We thus consider some conditions on the uniform
consistency of the density estimator $p_n$. These assumptions
are satisfied for the usual kernel or wavelets estimators of the
transition density.

\begin{enumerate}
\item[\textbf{H3}] For a sequence of nonnegative real numbers $(\alpha_{n})_{n\in\mathbb{N}}$ converging to $0$ as
$n\rightarrow\infty$, $p(x,y)$ is estimated by $p_{n}(x,y)$ at the rate $\alpha_{n}$
 for the mean square error when error is measured by the $L^{\infty}$ loss over $S\times S$:
\[
\mathbb{E}_{\nu}\left[\sup_{(x,y)\in S\times S}
|p_{n}(x,x')-p(x,x')|^{2}\right]=\mathcal{O}_\nu(\alpha_{n}),\text{ as }n\rightarrow\infty.
\]
\item[\textbf{H4}] The minorizing probability $\Phi$ is such that $\inf_{x\in S}\phi(x)>0$.
\item[\textbf{H5}] The densities $p$ and $p_n$ are bounded over $S^2$ and $\inf_{x,y\in S}p_n(x,y)/\phi(y)>0$.
\end{enumerate}

Since the choice of $\Phi$ is left to the statistician, we can use for instance the uniform distribution overs $S$,
even if it may not be optimal to do so. In such a case, \textbf{H4} is automatically
satisfied. Similarly, it is not difficult to construct an estimator
$p_n$ satisfying the constraints of \textbf{H5}.

Results of the previous section can then be extended to Harris chains:
\begin{theorem}\label{main}
Let $\mu$ be the invariant measure of the chain, and $\theta_0\in
\mathbb{R}^p$ be the parameter of interest, satisfying
$\mathbb{E}_\mu[m(X,\theta_0)]=0$.
Consider $A=S\times\{1\}$ an atom of the split chain,
$\tau_A$ the hitting time of $A$ and $B=(X_1,\cdots,X_{\tau_A})$. Assume
the hypotheses \textbf{H3}, \textbf{H4} and \textbf{H5}, and suppose
that $\mathbb{E}_A[M(B,\theta_0)M(B,\theta_0)']$ is of full rank.
\begin{enumerate}
\item[(a)]If \textbf{H0}$(S,4,\nu)$ and \textbf{H0}$(S,2)$ holds
as well as \textbf{H1}$(S,4,\nu,m)$ and \textbf{H1}$(S,2,m)$, then we
have in the just-identified case ($r=p$):
$$
2\hat r_{n}(\theta_0)\xrightarrow[n\to\infty]{\mathcal L}\chi_p^2
$$
and therefore
          $$
          \widehat{C}_{n,\alpha}=\left\{ \theta \in \mathbb{R}^{p}\left|
          \;2\cdot\hat{r}_{n}(\theta)\leq F^{-1}_{\chi_p^2}(1-\alpha)\right.\right\}.
          $$
is an asymptotic confidence region of level $1-\alpha$.
\item[(b)]
Under the additional assumptions~\textbf{H2}(a), \textbf{H2}(b) and
\textbf{H2}(d),
$$
\hat\theta=\arg\inf_{\theta\in\Theta} \{\hat r_n(\theta)\}
$$
is a consistent estimator of $\theta_0$.
If in addition~\textbf{H2}(c) holds, then
$\sqrt n (\hat\theta-\theta_0)$ is asymptotically normal.
\item[(c)]
In the case of over-identification ($r>p$),
we have:
$$
\widehat W_{1,n}(\theta_0)= 2\hat r_n(\theta_0)-2\hat
r_n(\hat\theta)\xrightarrow[n\to\infty]{\mathcal L}\chi_{p}^{2}
$$
and
$$
\widehat C^1_{n,\alpha}=\left\{ \theta \in \mathbb{R}^{p}
\left|\widehat W_{1,n}(\theta) \leq F^{-1}_{\chi_p^2}(1-\alpha
)\right.\right\},
$$
is an asymptotic confidence region of level $1-\alpha$. The moment
equation~(\ref{atomic_moment}) can be tested by using the following
convergence in law:
$$2\hat r_n(\hat\theta)\xrightarrow[n\to\infty]{\text{under~(\ref{atomic_moment})}}\chi^2_{r-p}.$$
\item[(d)]
Let $\theta=(\gamma,\beta)'$, where $\gamma\in\mathbb R^q$
and $\beta\in\mathbb R^{p-q}$. Under the hypotheses $\gamma=\gamma_0$,
$$
\widehat W_{2,n}(\gamma_0)=
2\inf_{\beta}\hat r_n((\gamma_0,\beta)' )-2\hat r_n(\hat\theta)
\xrightarrow[n\to\infty]{\mathcal L}\chi^2_{q}$$ and then
$$
\widehat C^2_{n,\alpha}=\left\{ \gamma \in
\mathbb{R}^{q}\left|\widehat W_{2,n}(\gamma) \leq
F^{-1}_{\chi_q^2}(1-\alpha )\right.\right\},
$$
is an asymptotic confidence region of level $1-\alpha$ for the
parameter of interest $\gamma$.
\end{enumerate}
\end{theorem}

\section{Some simulation results}\label{illu}

\subsection{Illustrative example}
We introduce here an example in order to illustrate the method in a very simple setting and compare ReBEL with BEL~\citep{kitamura97} in a situation that favor none.

We consider a AR(1), which is also a Markov chain of order 1, defined as follows:
$$
X_0=0,\ \varepsilon_i\sim \mathcal{U}([-\sqrt{12},\sqrt{12}])\text{ and } X_i=0.9 X_{i-1}+\varepsilon_i,
$$
$\varepsilon_t$ being i.i.d. uniformly distributed random variables of mean $0$ and variance $1$. Since $0.9<1$, the chain is recurrent.
The parameter of interest is the mean of the chain, $\mu=0$.

We choose in the simulations a small set $S=[-a;a]$. On each simulation, $a$ is chosen from a small grid in order to maximize the number of regeneration times as explained in section~\ref{algo}. ReBEL algorithm is then implemented and coverage probabilities are calculated for the nominative level 95\%.

We also compute BEL coverage probabilities, using non-overlapping blocks of constant length the integer part of $n^{1/3}$, where $n$ is the data set length.

We get the following results, for 10 000 replications:

$$
\text{table~\ref{toy_example_tab} should approximatively here}
$$

On these simulations, ReBEL seems to be better fitted. This may be due to the fact that ReBEL' small set length is data driven whereas BEL's blocks length is constant over the replications. In the following section, we set the small set once for all the replication.

\subsection{Estimation of the threshold crossing rate of a TGARCH}\label{tgarch}

The aim of this section is to show that ReBEL can be adapted to complex data and can outperform competing methods. Some applications of empirical likelihood to dependent data have been carried out, such as \cite{li-wang03} on Stanford Heart Transplant data or \cite{ow01} on bristlecone pine tree rings. In his book, Owen motivates his use of empirical likelihood to study the tree rings data set by its asymmetry: ``we could not capture such asymmetry in an AR model with normally distributed errors''~\citep[][page~168]{ow01}.

To motivate the use of empirical likelihood, we propose here to generate data sets with strong asymmetry properties to illustrate the applicability of the method.
For this, we consider a family of models introduced to study financial data, the TGARCH~\citep{zakoian93}. This model has been designed to handle non symmetric data, such as stock return series in presence of asymmetry in the volatility.
We think in particular to applications on modeling electricity prices series \citep{cornec-harari08}.
These series are very hard to model because of their very asymmetric behavior and  because of the presence
of very sharp peaks alternating with periods of low volatility. Application of ReBEL to these series seems to
be promising, see \cite{cornec-harari08}.

The data generating process is the following:
\begin{equation*}\left\{\begin{array}{rll}
X_i=&0.97 X_{i-1} +\varepsilon_i&\text{ with }X_0=0,\\
\varepsilon_i=&\sigma_i \nu_i&\text{ with }\nu_i\sim\text{NID}(0,1),\\
\sigma_i=&1+0.5 |\varepsilon_{i-1}|+0.4\varepsilon_{i-1}^+&\text{ with }\varepsilon_0=0.
\end{array}\right.\end{equation*}

where the $\nu_i$ are standard normal random values independent of all other random variables
and $x^+$ is the positive part of $x$: $x^+=\max\{0,x\}$. Of course, in the following, this generating mechanism is considered unknown.
Retrieving the underlying mechanism by just looking at the data is a difficult task and this motivates the use of
a non parametric approach in this context.

It is straightforward that $(X,\varepsilon)$ is a Markov chain of order 1. As $\varepsilon_{i-1}=X_{i-1}-0.97X_{i-2}$, it is immediate that $X$ is a Markov chain of order 2.
ReBEL algorithm can then be applied to $\mathbf X_i =(X_i,X_{i-1})$, which is a Markov chain of order 1.

In practice, the order $k$ of the Markov chain is unknown and is therefore to be estimated. We propose the following heuristic procedure to estimate the order:
\begin{itemize}
\item[(1)] Suppose $k=1$.
\item[(2)] Build the block according to Algorithm~\ref{algo2}.
\item[(3)] Evaluate the moment condition over the blocks: $Y_j=M(\widehat B_j,\theta)$.
\item[(4)] Perform a test of independence (or at least of non correlation) of the $(Y_1,\dots,Y_{l_n-1})$, for example by testing the nullity of $\rho$ given by $Y_i=\rho Y_{i-1}+\nu_i$. Other tests may be
considered as well, such as tests based on kernel estimators of the density.
\item[(5)] If the independence (or non correlation) is rejected, set $k=k+1$ and restart at point 2.
\end{itemize}

In order to apply~\cite{kitamura97}'s Block Empirical Likelihood (BEL), $X$ must be weakly dependent.
As the sum of the coefficients of $|\varepsilon_{i-1}|$ and $\varepsilon_{i-1}^+$ is smaller than $1$, the volatility of the data generating process is contracting. Therefore one can easily check the weak dependence of the process.

\subsection{Confidence intervals}
We are interested in estimating the probability of crossing a high threshold.
This is an interesting problem because of the asymmetry
of the data and a problem of practical interest for electricity prices. Indeed, production means are
only profitable above some level. The probability of crossing the profitability threshold is therefore
essential to estimate. The parameter of interest is defined here as:
$$
\theta_0=\mathbb E_\mu\left[\1_{\left\{ X_i\geq10\right\}}\right]=\mathbb P_\mu\left( X_i\geq10\right).
$$
and its value (estimated on a simulated data set of size $10^6$) is $\theta_0=0.1479$. A first advantage of ReBEL is that such a parameter, defined with respect to the underlying invariant measure $\mu$, is naturally handled by this method, whereas no unbiased estimating equation is available for BEL.

\begin{figure}
\includegraphics[width=16cm] {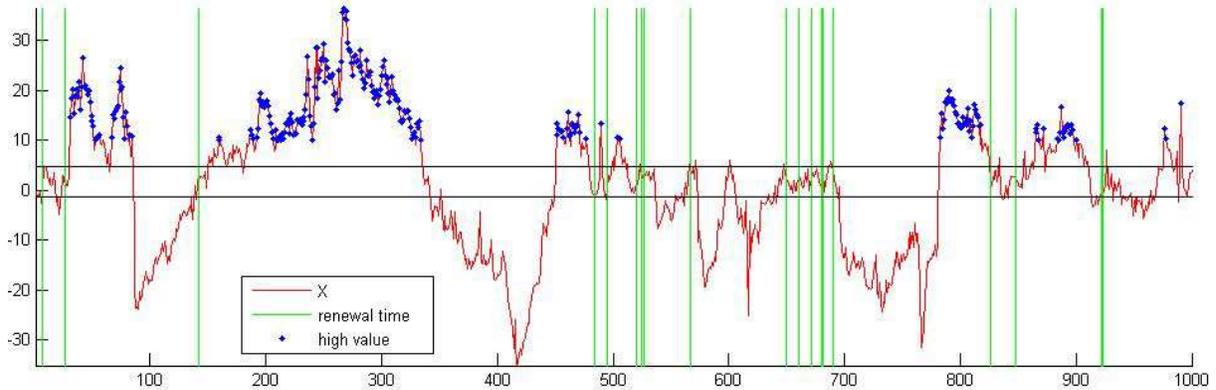}
\caption{The plain curve is a chain of length 1000. The horizontal lines limit the small set.
The 18 renewal times are marked by vertical lines. High values are marked by a dot.}\label{chaine}
\end{figure}

We simulate a data set of length 1000 and perform a test to estimate the order of the chain. We build an estimator $p_n$ of the transition density $p$ based on Gaussian kernels. The hypothesis $q=1$ is rejected whereas $q=2$ is not. As the chain is then considered 2-dimensional, we consider a small set of the form $S^2$ where $S$ is an interval. The interval $S$ has been chosen empirically to maximize the number of blocks and is equal to $[-1.3;4.7]$. It is set once and for all and is not updated at each replication. On the graphic corresponding to one simulation, $\mathbf X$ is in the small set $S^2$ when the trajectory of $X$ is in between the 2 plain black lines $y=-1.3$ and $y=4.7$ for two consecutive times.
For $i$ such that $\mathbf X_i$ visits $S^2$, we generate a Bernoulli $B_i$ as in Algorithm~\ref{algo2}, and if $B_i=1$, $i$ is a approximate renewal time. On the simulation, $S$ is visited 231 times, leading to 18 renewal times, marked by a vertical green line.
$$
\text{Figure~\ref{chaine} should be approximately here}
$$
The block length adapts to the local behavior of the chain: regions of low volatility lead to small blocks
(between 500 and 700) whereas regions with high values lead to larger blocks (like the 142-484 block). It can be noticed that high values concentrate in few blocks, because the dependence is well captured by Algorithm~\ref{algo2}.
BEL procedure leads to constant length block which cannot adapt to the dependence structure.
As suggested by~\cite{hall-horowitz-jing95}, the BEL blocks used in the following are of length $n^{1/3}=10$ and
then the chain is divided into 100 non overlapping blocks. The overlapping block perform poorly and won't be considered in the following.

Now that we have ReBEL approximately regenerative blocks, we can apply Theorem~\ref{main}(a) to obtain a confidence interval for $\theta$. We give a BEL confidence interval as well for comparison. We also consider two simpler methods as references for the performances of ReBEL: the simple sample mean 
$mean=n^{-1}\sum_i \1_{\left\{ X_i\geq10\right\}}$ and the mean over the regenerative blocks 
$$
trunc=\frac{\sum_{k=1}^{l_n}\sum_{i=\hat\tau(k)+1}^{\hat\tau{k+1}}\1_{\{X_i\geq10\}}}
{\sum_{k=1}^{l_n}(\hat\tau{(k+1)}-\hat\tau(k))}
=\frac{\sum_{k=1}^{l_n}\sum_{i=\hat\tau(k)+1}^{\hat\tau{k+1}}\1_{\{X_i\geq10\}}}
{\hat\tau{(l_n+1)}-\hat\tau(1))}.
$$
The simple mean do not deal with the dependence and we expect it to perform poorly. The second reference method $trunc$ uses the splitting technique in its expression, but in practice it only differs from $mean$ by the fact that it discards the first and last blocks.

An important point here is that to build confidence intervals with these two methods, an estimator of the variance is needed. In fact, if these estimators seem much simpler than BEL and ReBEL, the difficulty is mainly transferred to the estimation of their variances. This issue is difficult in a general dependence setting. In the applications, we used a bootstrap estimator of the variance of $mean$ and $trunc$ according to~\cite{GK96}.

\begin{figure}\begin{center}
\includegraphics[width=15cm] {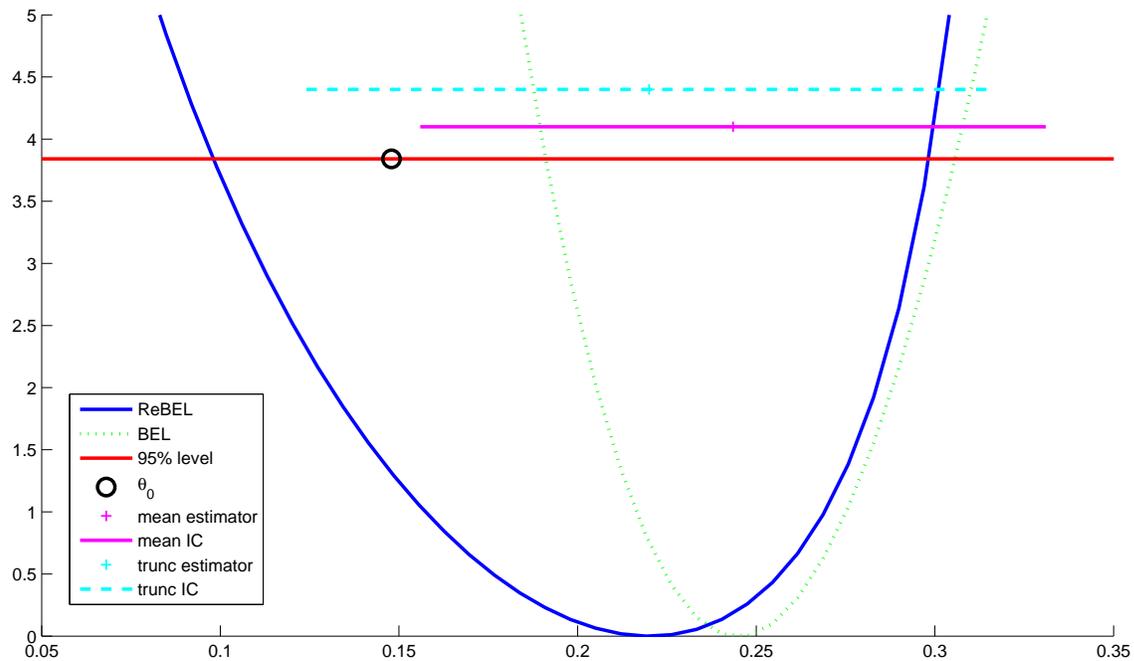}
\caption{The plain curve gives the ReBEL likelihood, whereas the dotted curve shows
the BEL likelihood. The red horizontal line marks the 95\% level and $\theta_0$ is marked by a circle on that line. The $mean$ CI is the magenta segment whereas the $trunc$ CI is the larger dotted cyan segment.}\label{IC}
\end{center}\end{figure}

Having in mind that difficulty, it is important to stress that ReBEL and BEL confidence intervals do not rely on an estimation of the variance of the estimator. This property is well-known for methods based on empirical likelihood that automatically estimate a variance at each point of the confidence interval, see for example the Continuously updated GMM~\cite{bonnal-renault07}. Additional results on the self-normalized properties of these methods have been investigated in~\cite{ber-gaut-hara08}.
$$
\text{Figure~\ref{IC} should be approximately here}
$$
$Mean$ and BEL estimators and confidence intervals appear biased to the right. This is most likely due to the effect of data from the first and last blocks, discarded by ReBEL and $trunc$.  It can be noticed that BEL blocks being more numerous, the confidence interval is tighter for BEL than for ReBEL.

To compare the considered methods, we also compute coverage probabilities and type-II errors
(which is equivalent to power in terms of test) of confidence intervals with nominal level 95\%.
To test the behavior under the alternative, we evaluate the statistics at the erroneous points $\theta=\theta_0+5/\sqrt n$ and $\theta=\theta_0+10/\sqrt n$ and check if the null hypotheses if rejected or not.

The 2 000 simulation results are summarized in Table~\ref{power}, for $n=$1000, 5000 and 10000.
$$
\text{Table~\ref{power} should be approximately here}
$$
Globally, ReBEL's coverage probabilities are better than BEL's, whereas its type-II error are bigger.
This is coherent with Figure~\ref{IC}: ReBEL confidence interval leads to better
coverage probabilities but is larger than BEL's (and therefore type-II errors are bigger for ReBEL).
$Mean$ and $trunc$ perform well for $n=1000$ but show some limits for $n=5000$ and $10000$.
It seems that ReBEL is the only method converging to the nominal level 95\%.

\begin{figure}\begin{center}\includegraphics[width=13cm] {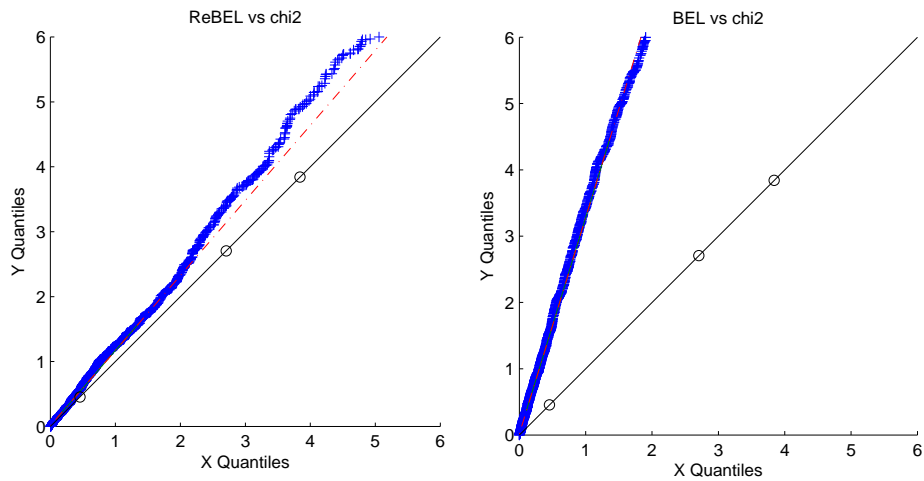}
\caption{QQ-plots of 10 000 Monte-Carlo repetitions of ReBEL statistic
versus $\chi^2_1$ quantiles. The solid reference line is the 45° line. The reference circles on that line mark
the 50\%, 90\% and 95\% levels. Data set length is $n=10\;000$.}\label{qqplot}
\end{center}\end{figure}

Coverage probabilities at other nominal level can also be investigated, and we make a Monte-Carlo experiment (10 000 repetitions) in order to confirm the adequacy to the asymptotic distribution achieved by the ReBEL algorithm. Data sets length are 10 000.
$$
\text{Figure~\ref{qqplot} should be approximately here}
$$
Figure~\ref{qqplot} shows the adequacy of the log likelihood to the asymptotic distribution given by Theorem~\ref{main}. The QQ-plots is almost linear and is close to the 45° line.

\section{Conclusion}

This paper propose an alternative point of view on dependent data sets and a corresponding semi-parametric methodology.
Random length blocks allow to adapt to the dependence structure of the data.
We have shown that ReBEL enjoys desirable properties corresponding to that of optimal reference methods for strong-mixing series. Simulations indicate that our algorithm at least competes with Kitamura's BEL when both methods can be applied.

This method seems to be a promising tool to handle dependent data when classical parametric models do not perform well, for example in presence of asymmetry and non normality of the innovations.


\bibliographystyle{gNST}
\bibliography{biblio}


\appendix\section{Proofs}
\subsection{Lemmas for the atomic case}

Denote $Y_j=M(B_j,\theta_0)$, $\overline{Y}=1/l_n\sum_{j=1}^{l_n}Y_j$ and define
$$
S^2_{l_n}=1/l_n\sum_{j=1}^{l_n}M(B_j,\theta_0)M(B_j,\theta_0)'=1/l_n\sum_{j=1}^{l_n}Y_jY_j'
 \text{ and }S^{-2}_{l_n}=(S^2_{l_n})^{-1}.
$$

To demonstrate Theorem~\ref{th_rebel}, we need 2 technical lemmas.
\begin{lemma}\label{lemma_s2}
Assume that $\mathbb{E}_A[M(B,\theta_0)M(B,\theta_0)']$
exists and is full-rank, with ordered eigenvalues $\sigma_p\geq \cdots\geq \sigma_1>0$.
Then, assuming \textbf{H0}$(1,\nu)$ and \textbf{H0}$(1)$, we have
$$S^2_{l_n}\to_{\nu}\mathbb{E}_A[M(B,\theta_0)M(B,\theta_0)'].$$
Therefore, for all $u\in\mathbb R^p$ with $\|u\|=1$,
$$
\sigma_1+o_\nu(1)\leq u'S^2_{l_n}u\leq \sigma_p+o_\nu(1).
$$
\end{lemma}
\paragraph{Proof:} The convergence of $S^2_{l_n}$ is a LLN for the sum of a random number of random variables,
and is a straightforward corollary of the Theorem~6 of~\citep[][chapter~5.2, page 131]{chow-teicher88}.

\begin{lemma}\label{lemma_max}
Assuming \textbf{H0}$(1,\nu)$, \textbf{H0}$(2)$ and \textbf{H1}$(2,m)$, we have
$$\max_{1\leq j\leq l_n} ||Y_j||=o_\nu(n^{1/2}).$$
\end{lemma}
\paragraph{Proof:} By \textbf{H1}$(2,m)$,
$$\mathbb{E}_A\left[\left(\sum_{i=1}^{\tau_1}\|m(X_i,\theta_0)\|\right)^2\right]<\infty,$$
and then,
$$\mathbb{E}_A[||Y_1||^2]=
\mathbb{E}_A\left[\left\|\sum_{i=1}^{\tau_1}m(X_i,\theta_0)\right\|^2\right]<\infty.$$
By Lemma A.1 of~\cite{bonnal-renault07},
the maximum of $n$ i.i.d. real-valued random variables with finite variance is $o(n^{1/2})$.
Let $Z_n$ be the maximum of $n$ independent copies of $||Y_1||$, $Z_n$ is then such as
$Z_n=o_\nu(n^{1/2})$.
As $l_n$ is smaller than $n$, $\max_{1\leq j\leq l_n} ||Y_j||$ is bounded by $Z_n$ and therefore,
$\max_{1\leq j\leq l_n} ||Y_j||=o_\nu(n^{1/2})$.

\subsection{Proof of Theorem~\ref{th_rebel}}
The likelihood ratio statistic $r_n(\theta_0)$ is the supremum
over $\lambda\in\mathbb{R}^p$ of $\sum_{j=1}^{l_n}\log(1+\lambda'Y_j)$.
The first order condition at the supremum $\lambda_n$ is
then:
\begin{equation}
1/l_n\sum_{j=1}^{l_n}\frac{Y_j}{1+\lambda_n' Y_j}=0.\label{1ord_rebel}
\end{equation}
Multiplying by $\lambda_n$ and using $1/(1+x)=1-x/(1+x)$, we have
$$
1/l_n\sum_{j=1}^{l_n}(\lambda_n' Y_j)\left(1-\frac{\lambda_n' Y_j}{1+\lambda_n' Y_j}\right)=0,
\text{ and then }
\lambda_n'\overline{Y}=1/l_n\sum_{j=1}^{l_n}\frac{\lambda_n' Y_j Y_j'\lambda_n}{1+\lambda_n' Y_j}.
$$

Now we may bound the denominators $1+\lambda_n' Y_j$ by $1+||\lambda_n|| \max_j ||Y_j||$ and then
$$
\lambda_n'\overline{Y}=1/l_n\sum_{j=1}^{l_n}\frac{\lambda_n' Y_j Y_j'\lambda_n}{1+\lambda_n' Y_j}
\geq\frac{\lambda_n'S^2_{l_n}\lambda_n}{(1+||\lambda_n|| \max_j ||Y_j||)}
$$
Multiply both sides by the denominator,
$\lambda_n'\overline{Y}(1+||\lambda_n|| \max_j ||Y_j||)\geq\lambda_n'S^2_{l_n}\lambda_n$ or
$$
\lambda_n'\overline{Y}\geq\lambda_n'S^2_{l_n}\lambda_n-||\lambda_n|| \max_j ||Y_j||\lambda_n'\overline{Y}.
$$
Dividing by $||\lambda_n||$ and setting $u=\lambda_n/||\lambda_n||$, we have
\begin{equation}
u'\overline{Y}\geq||\lambda_n||\left[u'S^2_{l_n}u- \max_j ||Y_j|| u'\overline{Y}\right].\label{ineq}
\end{equation}
Now we control the terms between the square brackets.
First, by Lemma~\ref{lemma_s2}, $u'S^2_{l_n}u$ is bounded between $\sigma_1+o_\nu(1)$
and $\sigma_p+o_\nu(1)$.
Second, by Lemma~\ref{lemma_max}, $\max_j ||Y_j||=o_\nu(n^{1/2})$.
Third, the CLT applied to the $Y_j$'s gives $\overline{Y}=\mathcal{O}_\nu(n^{-1/2})$.
Then, inequality~(\ref{ineq}) gives
$$
\mathcal{O}_\nu(n^{-1/2})\geq||\lambda_n||\left[u'S^2_{l_n}u-o_\nu(n^{1/2})\mathcal{O}_\nu(n^{-1/2})\right]=
||\lambda_n||(u'S^2_{l_n}u+o_\nu(1)),
$$
and $||\lambda_n||$ is then $\mathcal{O}_\nu(n^{-1/2})$.

Using the first order condition~(\ref{1ord_rebel}) as well as
the equality $1/(1+x)=1-x+x^2/(1+x)$, we get
$$
0=1/l_n\sum_{j=1}^{l_n} Y_j\left(1-\lambda_n' Y_j +\frac{(\lambda_n' Y_j)^2}{1+\lambda_n' Y_j}\right)
=\overline{Y}-S^2_{l_n}\lambda_n+1/l_n\sum_{j=1}^{l_n} \frac{Y_j(\lambda_n' Y_j)^2}{1+\lambda_n' Y_j}.
$$
The last term is $o_\nu(n^{-1/2})$ by Lemma~A.2 of~\cite{bonnal-renault07} and then
$\lambda_n=S^{-2}_{l_n}\overline{Y}+o_\nu(n^{-1/2})$.

Now, developing the $\log$ up to the second order,
$$
2r_n(\theta_0)=2\sum_{j=1}^{l_n}\log(1+\lambda_n'Y_j)
=2l_n\lambda_n'\overline{Y}-l_n\lambda_n'S^2_{l_n}\lambda_n+2\sum_{j=1}^{l_n}\eta_j,
$$
where the $\eta_i$ are such that, for some positive $B$ and with probability tending to 1,
$|\eta_j|\leq B|\lambda_n'Y_j|^3$.
Since, by Lemma~\ref{lemma_max}, $\max_j ||Y_j||=o_\nu(n^{1/2})$,
$$
\sum_{j=1}^{l_n}\|Y_j\|^3\leq n\max_j ||Y_j||\left(\frac1{l_n}\sum_{j=1}^{l_n}\|Y_j\|^2\right)
=n o_\nu(n^{1/2})\mathcal{O}_\nu(1)=o_\nu(n^{3/2})
$$
from which we find
$$
2\sum_{j=1}^{l_n}\eta_j\leq B\|\lambda_n\|^3\sum_{j=1}^{l_n}\|Y_j\|^3
=\mathcal{O}_\nu(n^{-3/2})o_\nu(n^{3/2})=o_\nu(1).
$$
Finally,
$$
2r_n(\theta_0)=2l_n\lambda_n'\overline{Y}-l_n\lambda_n'S^2_{l_n}\lambda_n+o_\nu(1)
=l_n\overline{Y}'S^{-2}_{l_n}\overline{Y}+o_\nu(1)\xrightarrow[n\to\infty]{\mathcal L}\chi^2_p.
$$
This concludes the proof of Theorem~\ref{th_rebel}.

\subsection{Proof of Theorem~\ref{th_power}}

We keep the notations of the previous subsection. Note that instead of
$\mathbb E_A[Y_j]=0$, we have
$$
\mathbb E_A[Y_j]=\delta\frac{\mathbb E_A[\tau_A]}{\sqrt n}\sim\delta\sqrt{\frac{\mathbb E_A[\tau_A]}{l_n}},
$$
because $l_n\sim n/\mathbb E_A[\tau_A]$.
The beginning of the proof is similar to that of Theorem~\ref{th_rebel}:
the misspecification is not significant at first order:
$\overline Y$ remains $\mathcal O_\nu (n^{-1/2})$.
We obtain:
$$
2r_n(\theta_0)=l_n\overline{Y}'S^{-2}_{l_n}\overline{Y}+o_\nu(1).
$$
$\overline{Y}-\delta\sqrt{\mathbb E_A[\tau_A]/l_n}$ is asymptotically Gaussian with variance
$\mathbb{E}_A[M(B,\theta_0)M(B,\theta_0)']$, which is the limit in probability of $S^{-2}_{l_n}$.
Therefore,
$$
2r_{n}(\theta_0)\xrightarrow[n\to\infty]{\mathcal L} \chi^{'2}_p(\delta'\Sigma^{-1}\delta).
$$

\subsection{Proof of Theorem~\ref{atomic_estimator}}
In order to prove Theorem~\ref{atomic_estimator},
we use a result established by \cite{qin94}.
\begin{lemma}[Qin \& Lawless, 1994]\label{th_qin}
Let $Z,Z_1,\cdots,Z_n\sim F$ be i.i.d. observations in $\mathbb R^d$ and consider
a function $g:\mathbb R^d\times\mathbb R^p\to\mathbb R^r$ such that $\mathbb E_F[g(Z,\theta_0)]=0$.
Suppose that the following hypotheses hold:
\begin{enumerate}
\item[(1)] $\mathbb E_F[g(Z,\theta_0)g'(Z,\theta_0)]$ is positive definite,
\item[(2)] $\partial g(z,\theta)/\partial\theta$ is continuous and bounded in norm
by an integrable function $G(z)$ in a neighborhood $V$ of $\theta_0$,
\item[(3)] $||g(z,\theta)||^3$ is bounded by $G(z)$ on $V$,
\item[(4)] the rank of $\mathbb E_F\left[\partial g(Z,\theta_0)/\partial\theta\right]$ is $p$,
\item[(5)] $\cfrac{\partial^2 g(z,\theta)}{\partial\theta\partial\theta' }$ is continuous and bounded by $G(z)$ on $V$.
\end{enumerate}
Then, the maximum empirical likelihood estimator $\tilde\theta_n$ is
a consistent estimator and
$\sqrt{n}(\tilde\theta_n-\theta_0)$ is asymptotically normal with
mean zero.
\end{lemma}

Set
$$
Z=B_1=\left(X_{\tau_A(1)+1},\cdots,X_{\tau_A(2)}\right)\in\bigcup_{n\in\mathbb N}\mathbb R^n
$$ and
$g(Z,\theta)=M(B_1,\theta)$. Expectation under $F$ is then replaced by $\mathbb{E}_A$.
Theorem~\ref{atomic_estimator} is a straightforward application of the Lemma~\ref{th_qin}
as soon as the assumptions hold.

By assumption, $\mathbb{E}_A[M(B,\theta_0)M(B,\theta_0)']$ is of full rank. This implies (1).

By~\textbf{H2}(a), there is a neighborhood $V$ of $\theta_0$ and a function $N$ such that,
for all $i$ between $\tau_A+1$ and $\tau_A(2)$,
$\partial m(X_i,\theta)/\partial\theta$ is continuous on $V$ and bounded in norm by $N(X_i)$.
$\partial M(B_1,\theta)/\partial\theta$ is then continuous as a sum of continuous functions and is
bounded for $\theta$ in $V$ by $L(B_1)=\sum_{i=\tau_A(1)+1}^{\tau_A(2)}N(X_i)$.
Since $N$ is such that
$\mathbb{E}_{\mu}\left[N(X)\right]<\infty$, we have by Kac's Theorem,
$$
\mathbb{E}_{A}\left[\sum_{i=\tau_A(1)+1}^{\tau_A(2)}N(X_i)\right]/\mathbb E_A[\tau_A]
=\mathbb{E}_{A}[L(B_1)]/\mathbb E_A[\tau_A]<\infty.$$
The bounding function $L(B_1)$ is then integrable.
This gives assumption~(2). Assumption~(5) is derived from~\textbf{H2}(c) by the same arguments.

By~\textbf{H2}(d),
$\|m(X_i,\theta)\|^3$ is bounded by $N(X_i)$ for $\theta$ in $V$, and then
$$
\|M(B_1,\theta)\|^3\leq\sum_{i=\tau_A(1)+1}^{\tau_A(2)}\|m(X_i,\theta)\|^3
\leq\sum_{i=\tau_A(1)+1}^{\tau_A(2)}N(X_i)=L(B_1).
$$
Thus, $\|M(B_1,\theta)\|^3$ is also bounded by $L(B_1)$ for $\theta$ in $V$, and hypotheses~(3) follows.

By Kac's Theorem,
$$
\mathbb E_A[\tau_A]^{-1}\mathbb E_A[\partial M(B_1,\theta_0)/\partial\theta]
=\mathbb E_\mu[\partial m(X_i,\theta_0)/\partial\theta],$$
which is supposed to be of full rank by~\textbf{H2}(b).
Thus $\mathbb E_A[\partial M(B_1,\theta_0)/\partial\theta]$ is of full rank and this gives assumption~(4).
This concludes the proof of Theorem~\ref{atomic_estimator}.

Under the same hypotheses, Theorem~2 and Corollaries~4 and~5 of~\cite{qin94} hold.
They give respectively our Theorems~\ref{atomic_overid}, \ref{atomic_test} and~\ref{atomic_subvector}.

\subsection{Proof of Theorem~\ref{main}}

Suppose that we know the real transition density $p$. The chain
can then be split with the Nummelin technique as above. We get an
atomic chain $\widetilde{X}$. Let's denote by $B_j$ the blocks
obtained from this chain. The Theorem~(\ref{th_rebel}) can then be
applied to $Y_j=M(B_j,\theta_0)$.

Unfortunately $p$ is unknown 
and then we can not use the $Y_j$. Instead, we have the vectors
$\widehat{Y}_j=M(\widehat{B}_j,\theta_0)$, built on approximatively
regenerative blocks. To prove the Theorem~\ref{main}, we essentially
need to control the difference between the two statistics
$\overline{Y}=\frac{1}{l_n}\sum_{j=1}^{l_n}Y_j$ and
$\widehat{Y}=\frac{1}{\hat{l}_n}\sum_{j=1}^{\hat{l}_n}\widehat Y_j$.
This can be done by using Lemmas~(5.2) and~(5.3) in~\cite{BC04c}: under \textbf{H0}$(S,4,\nu)$,
we get
\begin{equation}
\left|\frac{\hat{l}_n}{n}-\frac{l_n}{n}\right|
=\mathcal{O}_{\nu}(\alpha_n^{1/2})\label{control_l}
\end{equation}
and under \textbf{H1}$(S,4,\nu,m)$ and \textbf{H1}$(S,2,m)$,
\begin{equation*}
\left\|\frac{\hat l_n}{n}\widehat{Y}-\frac{l_n}{n}\overline{Y}\right\|=
\left\|\frac{1}{n}\sum_{j=1}^{\hat{l}_n}\widehat{Y}_j-
\frac{1}{n}\sum_{j=1}^{l_n}Y_j\right\|=
\mathcal{O}_{\nu}(n^{-1}\alpha_n^{1/2}).
\end{equation*}
With some straightforward calculus, we have
\begin{equation}
\left\|\widehat Y-\overline Y\right\|\leq
\frac{n}{\hat{l}_n}\left\|\frac{\hat l_n}{n}\widehat{Y}-\frac{l_n}{n}\overline{Y}\right\|
+\left|\frac{l_n}{\hat{l}_n}-1\right| \left\|\overline{Y}\right\|.\label{trans}
\end{equation}
Since
$$
\left|\vphantom\int l_n-n/\mathbb{E}_A[\tau_A]\right|\xrightarrow[n\to\infty]{a.s.} 0,
$$
equation~(\ref{control_l}) gives
$$
\frac{n}{\hat{l}_n}=\frac{n}{l_n}\left(1+\frac{n}{l_n}\frac{\hat{l}_n-l_n}{n}\right)^{-1}
=\mathcal{O}_{\nu}(\mathbb{E}_A[\tau_A])
\left(1+\mathcal{O}_{\nu}(\mathbb{E}_A[\tau_A])\mathcal{O}_{\nu}(\alpha_n^{1/2})\right)^{-1}
=\mathcal{O}_{\nu}(\mathbb{E}_A[\tau_A])
$$
and
$$
\left|\frac{l_n}{\hat{l}_n}-1\right|=\frac{n}{\hat{l}_n}\left|\frac{\hat{l}_n-l_n}{n}\right|
=\mathcal{O}_{\nu}(\mathbb{E}_A[\tau_A])\mathcal{O}_{\nu}(\alpha_n^{1/2})
=\mathcal{O}_{\nu}(\alpha_n^{1/2}).
$$
From this and equation~(\ref{trans}), we deduce:
\begin{equation}
\left\|\widehat Y-\overline Y\right\|\leq
\mathcal{O}_{\nu}(\mathbb{E}_A[\tau_A])\mathcal{O}_{\nu}(n^{-1}\alpha_n^{1/2})+
\mathcal{O}_{\nu}(\alpha_n^{1/2})\mathcal{O}_{\nu}(n^{-1/2})
=\mathcal{O}_{\nu}(\alpha_n^{1/2}n^{-1/2}).\label{controlY}
\end{equation}
Therefore
$$
n^{1/2}\widehat Y=n^{1/2}\overline Y+n^{1/2}\left(\overline Y-\widehat Y\right)=
n^{1/2}\overline Y+\mathcal{O}_{\nu}(\alpha_n^{1/2}).
$$
Using this and the CLT for the $Y_i$, we show that $n^{1/2}\widehat Y$ is asymptotically Gaussian.

The same kind of arguments give a control on the difference between
empirical variances.
Consider
$$
\widehat S^2_{\hat l_n}=\sum_{j=1}^{\hat l_n}\widehat Y_j\widehat Y_j'
\text{ and }\widehat S^{-2}_{\hat l_n}=(\widehat S^2_{\hat l_n})^{-1}.
$$
By Lemma~(5.3) of~\cite{BC04c} we have, under
\textbf{H1}$(S,4,\nu,m)$ and \textbf{H1}$(S,2,m)$, $ \left\|\frac{\hat
l_n}{n}\widehat S^2_{\hat l_n}-\frac{l_n}{n}S^2_{l_n}\right\|=
\mathcal{O}_{\nu}(\alpha_n) $,
and then
\begin{equation}
\left\|\widehat S^2_{\hat l_n}-S^2_{l_n}\right\|
\leq
\frac{n}{\hat{l}_n}\left\|\frac{\hat l_n}{n}\widehat S^2_{\hat l_n}-\frac{l_n}{n}S^2_{ l_n}\right\|
+\left|\frac{l_n}{\hat{l}_n}-1\right| \left\|S^2_{l_n}\right\|
=\mathcal{O}_{\nu}(\alpha_n)+\mathcal{O}_{\nu}(\alpha_n^{1/2})
=o_{\nu}(1).\label{controlS}
\end{equation}

The proof of Theorem~(\ref{th_rebel}) is then also valid for the
approximated blocks $\widehat{B}_j$ and reduce to the study of the
square of a self-normalized sum based on the pseudo-blocks. We have
$\hat r_{n}(\theta_0)=\sup_{\lambda\in\mathbb{R}^p } \left\{
\sum_{j=1}^{\hat l_n}\log \left[1+\lambda'\widehat Y_j\right]
\right\}$. Let $\hat\lambda_n= -\widehat S^{-2}_{\hat
l_n}\widehat{Y}+o_{\nu}(n^{-1/2})$ be the optimum
value of $\lambda$, we have
$$
2\hat{r}_{n}(\theta_0)
=-2\hat l_n \hat\lambda_n'\widehat{Y}-\sum_{j=1}^{\hat l_n}(\hat\lambda_n' \widehat Y_j)^2+o_{\nu}(1)
=\hat l_n\widehat{Y}'\widehat S^{-2}_{\hat l_n}\widehat{Y}+o_{\nu}(1).
$$
Using the controls given by equations~(\ref{controlY}) and~(\ref{controlS}), we get
$$
2\hat{r}_{n}(\theta_0)=[l_n+\mathcal{O}_{\nu}(n\alpha_n^{1/2})]
\times\left[\overline{Y}'+\mathcal{O}_{\nu}\left(\sqrt\frac{\alpha_n}n\right)\right]
\times[S^{-2}_{l_n}+o_{\nu}(1)]
\times\left[\overline{Y}+\mathcal{O}_{\nu}\left(\sqrt\frac{\alpha_n}n\right)\right]
+o_{\nu}(1).
$$
Developing this product, the main term is $l_n  \overline{Y}
S^{-2}_{l_n} \overline{Y}\sim_{\nu}2r_{n}(\theta_0)$ and
all other terms are $o_{\nu}(1)$, yielding
$$
2\hat{r}_{n}(\theta_0)=l_n \overline{Y} S^{-2}_{l_n} \overline{Y}+o_{\nu}(1)
\xrightarrow[n\to\infty]{\mathcal L}\chi^2_{p}.
$$

Results~(b), (c) and~(d) can be derived from the atomic case
by using the same arguments.


\newpage

\begin{table}\label{toy_example_tab}\centering
\begin{tabular}{|c||c|c|}
\hline
$n$ & ReBEL & BEL\\
\hline
\hline
$250$&0.92&0.82\\
\hline
$500$&0.94&0.88\\
\hline
$1000$&0.94&0.91\\
\hline
\end{tabular}
\caption{Cover rates of confidence intervals for the mean of an AR(1). Comparison of ReBEL and BEL for different data set lengths. Nominal level is 0.95\;.}
\end{table}

\newpage

\newpage

\newpage
\begin{table}
$$
\begin{array}{|c||c|c|c|c||c|c|c|c||c|c|c|c|}
\hline
&\multicolumn{4}{c||}{\theta=\theta_0}&\multicolumn{4}{c||}{\theta=\theta_0+5/\sqrt n}&\multicolumn{4}{c|}{\theta=\theta_0+10/\sqrt n}\\
\hline
\vphantom{\int^1} n & \text{ReBEL} & \text{BEL} & mean & trunc & \text{ReBEL} & \text{BEL} & mean & trunc
& \text{ReBEL} & \text{BEL}& mean & trunc\\
\hline\hline 
 1 000 & 54 & 55 & 58 & 58 & 24 & 13 & 14 & 20 & 07 & 02 & 02 & 06 \\
\hline
5 000 & 88 & 67 & 74 & 74 & 52 & 27 & 30 & 29 & 12 & 01 & 01 & 02 \\
\hline
10 000 & 92 & 70 & 76 & 77 & 59 & 31 & 34 & 33 & 11 & 02 & 02 & 02\\
\hline
\end{array}$$
\caption{Coverage probabilities and type-II errors (percent) under the null and two alternatives, for ReBEL against BEL, and 2 reference methods. Nominal level is 95\%}
\label{power}
\end{table}

\end{document}